\documentclass[ 12pt,reqno]{article}

\usepackage[usenames]{color}
\usepackage{graphicx}
\usepackage{amscd}
\usepackage[english]{babel}
\usepackage{color}
\usepackage{fullpage}
\usepackage{psfig}
\usepackage{graphics,amsmath,amssymb}
\usepackage{amsthm}
\usepackage{amsfonts}
\usepackage{latexsym}
\usepackage{epsf}
\usepackage{float}
\usepackage{MnSymbol}
\usepackage{framed} 
\usepackage{tikz} 

\usepackage[colorlinks=true,
linkcolor=webgreen,
filecolor=webbrown,
citecolor=webgreen]{hyperref}

\definecolor{webgreen}{rgb}{0,.5,0}
\definecolor{webbrown}{rgb}{.6,0,0}

\newcommand{\Hbk}{{\em Handbook}}  

\theoremstyle{plain}
\numberwithin{equation}{section}

 \newcommand{\seqnum}[1]{\href{https://oeis.org/#1}{\underline{#1}}}

\newcommand{\beql}[1]{\begin{equation}\label{#1}}
\newcommand{\eeq}{\end{equation}}

\begin{document}


\setcounter{page}{1}


\begin{center}

{\large\bf ``A Handbook of Integer Sequences'' Fifty Years Later } \\
\vspace*{+.2in}

N. J. A. Sloane, \\
The OEIS Foundation Inc.,\\
11 So. Adelaide Ave.,
Highland Park, NJ 08904, USA \\
Email:  \href{mailto:njasloane@gmail.com}{\tt njasloane@gmail.com}

\vspace*{+.1in}

January 16, 2023

\vspace*{+.1in}

\begin{abstract}
Until 1973 there was no database of integer sequences.  Someone coming across
 the sequence $1, 2, 4, 9, 21, 51, 127, \ldots$ would have had no way of
 discovering that it had been studied since 1870 (today these are
 called the Motzkin numbers, and form entry A001006 in the database).
 Everything changed in 1973 with the publication of {\em A Handbook of Integer Sequences},
 which listed $2372$ entries. This report describes the fifty-year evolution of the database from the
  \Hbk~to its present form as {\em The On-Line Encyclopedia of Integer Sequences} (or {\em OEIS}),
 which contains $360,000$ entries, receives a million visits a day, and has been cited
 $10,000$ times,  often with a comment saying  ``discovered thanks to the OEIS''.
\end{abstract}
\end{center}


\section{Introduction}\label{Sec1}
\begin{figure}[htb]
        \begin{minipage}[b]{0.30\linewidth} 
        \centering
\begin{tikzpicture}[scale=1.];
\draw [fill=black] (0,0) circle[radius=0.1];
\draw [fill=black] (-.1,-.1) rectangle (.1,.1);
\draw [fill=black] (0,1) circle[radius=0.08];
\draw [fill=black] (-1,2) circle[radius=0.08];
\draw [fill=black] (0,2) circle[radius=0.08];
\draw [fill=black] (1,2) circle[radius=0.08];
\draw [fill=black] (-1,3) circle[radius=0.08];
\draw [fill=black] (1,3) circle[radius=0.08];
\draw [thick] (0,0)--(0,2);
\draw [thick] (0,1)--(-1,2);
\draw [thick] (0,1)--(0,2);
\draw [thick] (0,1)--(1,2);
\draw [thick] (0,2)--(-1,3);
\draw [thick] (0,2)--( 1,3);
\end{tikzpicture}
        \end{minipage}
        \hspace{0.5cm}
         \begin{minipage}[b]{0.30\linewidth} 
        \centering
\begin{tikzpicture}[scale=.8];
\draw [thick] (0,0) circle[radius=2];
\draw [thick] (1.732, 1) -- (-1.732, -1);
\draw [thick] (1.732, -1) -- (-1.732, 1);
\draw [thick] (1.414, 1.414) -- (-.5, -1.936);
\draw [thick] (-1.414, 1.414) -- (.5, -1.936);
\end{tikzpicture}
        \end{minipage}
        \hspace{0.5cm}
       \begin{minipage}[b]{0.30\linewidth}
        \centering
        \includegraphics[width=1.0\textwidth]{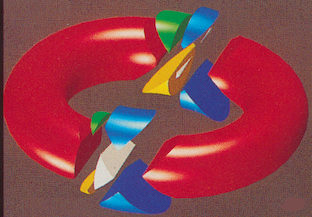}
        \end{minipage}
        \caption{Left: one of $48$ unlabeled rooted trees with $7$ nodes (the root node is at the bottom); center: four cuts of a pancake can produce $11$ pieces; right: three cuts of a bagel can produce $13$ pieces. }\label{Fig1}
        \end{figure}

Number sequences arise in all branches of science:
for example, $1, 1, 2, 4, 9, 20, 48, 115, \ldots$ gives the number
of rooted trees with
$n$ nodes (\seqnum{A000081},\footnote{Six-digit numbers prefixed by A refer to entries in the 
current version of the \Hbk,  {\em The On-Line Encyclopedia of Integer Sequences} \cite{OEIS}.}
see also Fig.~\ref{Fig1}),
and in daily life: how many pieces can you cut a pancake into with $n$
knife-cuts?  (The pieces need not all be the same size.) That one is easy: $1, 2, 4, 7, 11, 16, \ldots$, $n(n+1)/2 + 1$ (\seqnum{A000124}).
But what is the answer for cutting up an (ideal) bagel or doughnut? That is a lot harder: with a sharp knife you might
get a few terms, perhaps  $1, 2, 6, 13, \ldots$,  but probably not enough to guess the formula, which is
$n(n^2+3n+8)/6$ for $n>0$. For that
you would need to to consult the database: go to \href{https://oeis.org}{\tt https://oeis.org} and enter ``cutting bagel'',
or go directly to \seqnum{A003600}.

My fascination with these sequences  began in 1964 when I was a graduate student at Cornell University in Ithaca, NY, studying neural networks. I had encountered a sequence of 
numbers, $1, 8, 78, 944, 13800, \ldots$,
and I badly needed a formula for the $n$-th term, in order to determine 
the rate of growth of the terms (this
would indicate how long the activity in this very simple neural network would persist).
I will say more about this sequence in Section~\ref{SecHIS}.

I noticed that although several books in the Cornell library contained sequences 
somewhat similar to mine, as far as I could tell this particular sequence was not mentioned. 
I expected to have to analyze many related sequences, so in order to keep track of the 
sequences in these books, I started recording them on  $3" \times\, 5"$  file cards.

\begin{figure}[!ht]
   \centerline{\includegraphics[ width=1.5in]{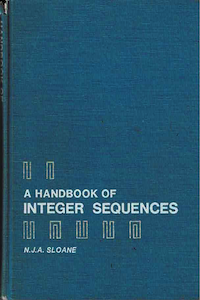}}
        \caption{Front cover of the \Hbk. The embossed figures show side views of
  the two ways of folding a strip
of three (blank) stamps, and the five ways of folding a strip of four stamps.
The full sequence begins $1, 1, 2, 5, 14, 38, 120, 353, 1148, 3527, \ldots$,
\seqnum{A001011}. No formula is known.}
        \label{FigHIS}
        \end{figure}

 The collection grew rapidly as I searched though more books, and once the word got out,
 people started sending me sequences.
 Richard Guy was an enthusiastic supporter right from the start.   In 1973 I formalized the collection as
 {\em A Handbook of Integer Sequences},  which was published by 
 Academic Press (Fig.~\ref{FigHIS}). It contained $2372$ entries.

Once the book appeared, the flood of correspondence increased, and it took twenty years to prepare
the next version.  Simon Plouffe helped a great deal, and in 1995 Academic Press published
our sequel, {\em The Encyclopedia of Integer Sequences}, with $5487$ entries.
From this point on the collection grew even more rapidly. I waited a year, until it had doubled in
size, and then put it on the Internet, calling it {\em The On-Line Encyclopedia of Integer Sequences}.

In the rest of this article I will first say more about the evolution of the database: the \Hbk~(\S\ref{SecHIS}),
the 1995 {\em Encyclopedia} (\S\ref{SecEIS}), the {\em On-Line Encyclopedia} (\S\ref{SecOEIS}), 
and the {\em OEIS Foundation} (\S\ref{SecOEISF}).
The next  sections describe the database itself:
what sequences are---or are not---included (\S\ref{SecWhat}),
how the database is used (\S\ref{SecUsed}),
the layout of a typical entry (\S\ref{SecEntry}), 
the arrangement of the entries (\S\ref{SecArrange}),
and a Fact Sheet (\S\ref{SecFacts}).
The final  sections describe some especially interesting sequences: 
Recam\'{a}n's sequence (\S\ref{SecRecaman}),
Iteration of number-theoretic functions (\S\ref{SecFunWithDigits}),
Gijswijt's sequence (\S\ref{SecGijs}),
Lexicographically Earliest Sequences (\S\ref{SecLES}),
The Stepping Stones problem (\S\ref{SecStones}),
Stained glass windows (\S\ref{SecScott}),
and Other sequences I would have liked to include (\S\ref{SecSpace}).

Several open questions are mentioned to which I would very much like to know the answers.

Notation. $a(n)$ denotes the $n$-th term of the sequence under discussion.
$\sigma(n)$ is the sum of the divisors of $n$ (\seqnum{A000203}).

%
%

\section{Evolution of the database}\label{SecHist}
\subsection{The {\em Handbook of Integer Sequences}}\label{SecHIS}
Once the collection had grown to a few hundred entries, I entered them
on punched cards,\footnote{These were never called ``punch cards'' (sic).
To anyone who worked with them in the 1960s, ``punch cards'' sounds like ``grill cheese'' (sic)
for ``grilled cheese'',  or ``barb wire'' (sic) for ``barbed wire'', both of which I have recently seen in print.} which
made it easier to check and sort them. The \Hbk~was typeset directly from the punched cards.
There were a few errors in the book, but almost all of them were caused by errors in the original publications.
Accuracy was a primary concern in that book, as it is today in the OEIS.

The book was an instant success.  It was, I believe,
the world's first dictionary of integer sequences (and my original title said
{\em Dictionary} rather than {\em Handbook}). Many people said ``What a great idea'', 
 and wondered why no one had done it before. 
 Martin Gardner recommended it in the {\em Scientific American} of July 1974.
Lynn A. Steen, writing in the {\em American Mathematical Monthly} said
``Incomparable, eccentric, yet very useful.
Contains thousands of `well-defined and interesting' infinite integer sequences 
together with references for each ...
If you ever wondered what comes after $1, 2, 4, 8, 17, 35, 71, \ldots$, this is the place to look it up''.

 Harvey J. Hindin, writing from New York City, exuberantly concluded a letter to me by saying: 
 ``There's the {\em Old Testament}, the {\em New Testament}, and the {\em Handbook of Integer Sequences}.''

I never did find the sequence that started it all in the literature,
but I learned P\'{o}lya's theory of counting, and with John Riordan's help found the answer,
which appears in \cite{RiSl69} and \seqnum{A000435}.

\subsection{The {\em Encyclopedia of Integer Sequences}}\label{SecEIS}
Following the publication of the \Hbk, a large amount of correspondence ensued, 
with suggestions for further sequences and updates to the entries.
By the early 1990's over a cubic meter of new material had accumulated. 
A Canadian mathematician, Simon Plouffe, offered to help in preparing a revised edition of the 
book, and in 1995 {\em The Encyclopedia of Integer Sequence}, by me and Simon Plouffe, was 
published by Academic Press. It contained $5487$ sequences, occupying $587$ pages.  
By now punched cards were obsolete, and the entries were stored on
magnetic tape.

\subsection{The {\em On-Line Encyclopedia of Integer Sequences}}\label{SecOEIS}
Again, once the book appeared, many further sequences and updates were 
submitted from people all over the world. I waited a year, until the size of the collection had 
doubled, to $10000$ entries, and then in 1996 I 
launched {\em The On-Line Encyclopedia of Integer Sequences}
 (now usually called simply the {\em OEIS}) on the Internet. From 1996 until October 26, 2009, 
it was part of my homepage on the AT\&T Labs website.

Incidentally, in 2004 the database was mentioned by the Internet website {\em slashdot}
 (``News for Nerds. Stuff that Matters"), and this brought so much traffic to my Bell Labs 
 homepage that it briefly crashed the whole Bell Labs website. My boss was quite proud of this, 
 since it was a rare accomplishment for the Mathematics and Statistics Research Center.

\subsection{{\em The OEIS Foundation}}\label{SecOEISF}
In 2009, in order to ensure the long-term future of the database,
 I set up a non-profit foundation,
{\em The OEIS Foundation Inc.}, a 501(c)(3) Public Charity, whose purpose is to own, 
maintain and raise funds to support {\em The On-Line Encyclopedia of Integer Sequences} 
or  {\em OEIS}.

On October 26, 2009, I transferred the intellectual property of {\em The On-Line Encyclopedia of Integer
 Sequences}  to the Foundation. 
 A new  OEIS with multiple editors  was launched on November 11, 2010.

Since then it has been possible for anyone in the world to propose a new sequence or an 
update to an existing sequence. To do this, users must first register, and then submissions are 
reviewed by  the editors  before they become a permanent part of the OEIS.
Technically the OEIS is now a ``moderated wiki''.

I started writing this article on November 11, 2022, noting that this marked
twelve years of successful operation of the online OEIS, and also that the database is in its $59$th year of existence.

%
%

\section{The database today}\label{SecToday}

\subsection{What sequences are included?}\label{SecWhat}

 From the very beginning the goal of the database has been to 
 include all ``interesting'' sequences of integers.
This is a vague definition, but some further examples may make it clearer.
The database  includes a huge number of familiar and unfamiliar sequences from mathematics
(the prime numbers $2, 3, 5, 7, 11, 13, \ldots$, \seqnum{A000040};
$60, 168, 360, 504, 660, 1092, \ldots$, the orders of noncyclic  simple groups, \seqnum{A001034}),
computer science ($0, 1, 3, 5, 8, 11, 14, \ldots$, the number of comparisons needed for 
merge sort, \seqnum{A001855}),
physics (see ``self-avoiding walks on lattices'', Ising model, etc., e.g. \seqnum{A002921}),
chemistry (the enumeration of chemical compounds was one of the 
motivations behind P\'{o}lya's theory of counting, see e.g. \seqnum{A000602}),
and not least, from puzzles and I.Q. tests  ($1, 8, 11, 69, 99, 96, 111, \ldots$, 
the ``strobogrammatic'' numbers, guess!, or see \seqnum{A000787};
$4, 14, 23, 34, 42, 50, 59, \ldots$, the numbered stops on the New York City A 
train subway, \seqnum{A011554}. 
That entry has links to a map and the train schedule).

Sequences that have arisen in the course of someone's work---especially if published---have always 
been welcomed.
On the other hand, sequences that have been proposed simply because they were missing 
from the database are less likely to be accepted.

There are a few hard and fast rules.
The sequence must be well-defined and the terms must not be time-dependent---if the next term 
is only  known to be either $14$ or $15$, for instance, then the sequence must end with the last term 
that is known for certain.  The sequence may not have any missing terms or gaps. 
In the case of Mersenne primes, for instance (\seqnum{A000043}) it is common for later primes
 to be known before all intermediate numbers having been tested. 
 The later primes get mentioned in comments, but they are not as part of the main sequence 
 until their position has been confirmed.
 
Very short sequences and sequences that are subsequences of many other sequences are not
accepted.  A sequence for which the only known terms are $2, 3, 5, 7$  would not be accepted 
since it is matched by a large number of existing sequences.
The definition may not involve an arbitrary but large parameter (primes ending in $1$ are fine,
 \seqnum{A030430}, but not primes ending in $2023$).
 
The OEIS \href{https://oeis.org/wiki/}{Wiki}
has a section listing additional  examples of what not to submit,
as well as a great deal of information
about the database that I won't repeat here, such as the meaning of the various keywords,
the definition of the ``offset'' of a sequence,  descriptions of the submission and editorial processes, 
and a list of over $10,000$ citations of the OEIS in the scientific literature.

Most OEIS entries give an ordered list of integers. But triangles of numbers are included 
by reading them row-by-row. For example, Pascal's triangle becomes $1, ~1,1,~$ $ 1, 2, 1, ~1, 3, 3, 1, \ldots$, 
\seqnum{A007318}.
Doubly-infinite square arrays are included by reading them by antidiagonals: 
the standard multiplication table for positive integers becomes $1,~ 2,2,~ $ $3, 4, 3, ~$ $4, 6, 6, 4, \ldots$, \seqnum{A003991}.

Sequences of fractions are included as a linked pair giving the numerators and denominators 
separately (the Bernoulli numbers are \seqnum{A027641}/\seqnum{A027642}).
Important individual  real numbers are included by giving their decimal or continued fraction 
expansions (for $\pi$
see \seqnum{A000796} and \seqnum{A001203}).
A relatively small number of sequences of nonintegral real numbers are included by 
rounding them to the nearest integer, or by taking floors or 
ceilings (the  imaginary parts of the zeros of Riemann's zeta function give \seqnum{A002410}).

Two less obvious sources for sequences are binomial coefficient identities and number-theoretic inequalities. 
The values of either side of the identity
 $$\sum_{k=0}^{n} \binom{2n}{k}^2 ~=~ \frac{1}{2} \binom{4n}{2n} - \frac{1}{2} \binom{2n}{n}^2$$
 \cite[(3.68)]{Gou72} give \seqnum{A036910}.
 From the inequality  $\sigma(n) <n \sqrt{n}$ for $n>2$,
 \cite[Sect.~III.1.1.b]{MSC96}, we get
 the integer sequence  $\lfloor n \sqrt{n} \rfloor - \sigma(n)$, \seqnum{A055682}.
 The point is that if you want to know if this inequality is known, you look up the 
 difference sequence, 
 and find \seqnum{A055682}  and a reference to the proof.  
 Many more sequences of these two types should be added to the database.

\subsection{How the database is used}\label{SecUsed}

The main applications of the database are in identifying sequences or in finding out the current status of a known sequence. Barry Cipra has called it a mathematical analogue of a ``fingerprint file''.
You encounter a number sequence, and wish to know if anyone has ever come across it before. 
 
If your sequence is in the database, the reply will give a definition, the first $50$
or so terms, and,
when available, formulas, references, computer code for producing the sequence, 
links to any relevant web sites, and so on.  

\begin{figure}[!ht]
\centerline{\includegraphics[ width=6in]{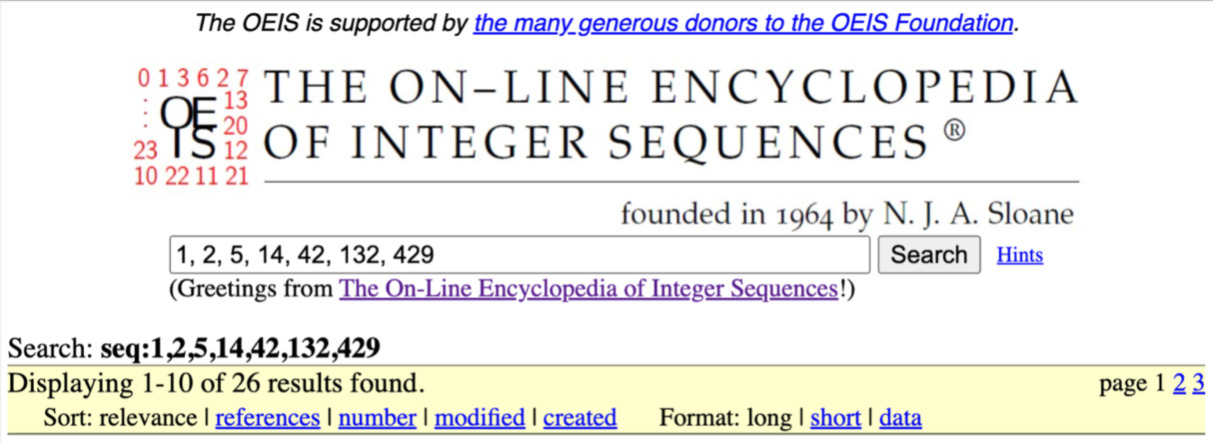}}
\caption{The result of submitting  
$1, 2, 5, 14, 42, 132, 429$ to the database. This figure shows the banner at the top of
the reply. There are $26$ matches, ranked in order of importance, the top match being the one we want, 
the  Catalan numbers. A shortened version of the top match is shown in the next figure.}
 \label{FigA108a}
\end{figure}

Figures~\ref{FigA108a} and \ref{TableA108b} show what happens if you submit 
$1, 2, 5, 14, 42, 132, 429$,  the first few Catalan numbers,
 one of the most famous sequences of all.

\begin{figure}[!ht]
\footnotesize
\begin{tabular}{|l|l|} \hline
A000108 & Catalan numbers: $C(n) = {\rm binomial}(2n,n)/(n+1) = (2n)!/(n!(n+1)!).$ \\ 
                & (Formerly M1459 N0577) \\  
 DATA      & 1, 1, 2, 5, 14, 42, 132, 429, 1430, 4862, 16796, 58786, ... \\
 COMMENTS & These were formerly sometimes called Segner numbers. \\
 & A very large number of combinatorial interpretations are known - \\
 & ~~see references,  esp. R. P. Stanley, Catalan Numbers, Camb., 2015.\\
  & This is probably the longest entry in the OEIS, and rightly so.\\
  & The solution to Schr\"{o}der's first problem: number of ways to insert $n$ pairs \\
  & ~~of parentheses in a word of $n+1$ letters. E.g., for $n=2$ there are 2 ways:  \\
  &  ~~((ab)c) or (a(bc)); for n=3 there are 5 ways: ((ab)(cd)), (((ab)c)d), ... \\
  & ... \\
  REFERENCES 
  & The large number of references and links demonstrates the ubiquity \\
  & ~~of the Catalan numbers. \\
& R. Alter, Some remarks and results on Catalan numbers, pp. 109-132 \\
& ~~ in Proc. Louisiana Conf.  Combinatorics, Graph Theory and \\
& ~~Computer Science. Vol. 2, edited R. C. Mullin et al., 1971.\\
 & M. Bona, ed., Handbook of Enumerative Combinatorics, CRC Press, 2015 \\
& L. Comtet, Advanced Combinatorics, Reidel, 1974, p. 53. \\
& J. H. Conway \& R. K. Guy, The Book of Numbers, Springer, 1995,  96-106. \\
& ... \\
LINKS & Robert G. Wilson v, Table of n, a(n) for n = 0..1000 \\
& ... \\
& F. R. Bernhart, Catalan, Motzkin and Riordan numbers, Disc. Math., \\
& ~~Vol. 204, No. 1-3 (1999), 73-112.\\
& ... \\
& W. G. Brown, Historical Note on a Recurrent Combinatorial Problem, \\
& ~~ Amer. Math. Monthly,. 72, No. 9 (1965), 973-977.\\
& ... \\
FORMULA & Recurrence: $a(n) = 2*(2*n-1)*a(n-1)/(n+1)$ with $a(0) = 1$.\\
& ... \\
MAPLE & A000108 := n->binomial(2*n, n)/(n+1); \\
& ... \\
MATHEMATICA & A000108[n\_] := (2 n)!/n!/(n+1)! \\
& ... \\
PARI & a(n)=binomial(2*n, n)/(n+1) \\
& ... \\
KEYWORD & core,nonn,easy,nice \\
AUTHOR & N. J. A. Sloane \\
 \hline
\end{tabular}
\caption{The entry for the Catalan number \seqnum{A000108}. The full entry has over 750 lines, 
which have been  edited  here to show samples of the different fields.}
\label{TableA108b}
\normalsize
\end{figure}

I could have chosen a simpler example, like the Fibonacci numbers, but I have a particular reason for 
choosing the Catalan numbers.  When the OEIS was new, people would sometimes say to me that they had a sequence 
they were trying to understand, and would I show them how to use the database. 
At least twice when I used the Catalan sequence as an illustration, they said, why, 
that is my sequence, how on earth did you know? It was no mind-reading trick, 
the Catalan numbers are certainly the most common sequence that people don't know about.
This entry is  the longest---and one of the most important---in the whole database.

If we do not find your sequence in the database, we will send you a message inviting 
you to submit it (if you consider it is of general interest), so that the next person who comes 
across it will be helped, and your name will go on record
as the person who submitted it.

The second main use of the database is to find out the latest information about a particular sequence.

Of course we cannot hope to keep all $360000$ entries up-to-date.   But when a new paper is published that mentions the OEIS, Google will tell us, and we then add links to that paper from any sequence that it mentions.
People have told us that this is one of the main ways they use the OEIS.  After all, even a specialist in (say) permutation groups cannot keep track of all the papers published worldwide in that area.
And if a paper in a physics journal happens to mention a number-theoretic sequence, for example, that is unlikely to be noticed by mathematicians.

There are also many other ways in which the database has proved useful.  

For example, it is an excellent source of problems to work on.  The database is  constantly being updated.
Every day we get  thirty to fifty submissions of new sequences, and an equal number of
comments on existing entries (new formulas, references, additional terms, etc.).
The new sequences  are often sent in by non-mathematicians,
and are a great source of problems. You can see the current submissions at \url{https://oeis.org/draft}.
Often enough you will see a sequence that is so interesting you want to drop everything  and work on it.
And remember  that we are always in  need of more volunteer editors.  In fact
 anyone who has registered with the OEIS can suggest edits, you do not even need to be an official editor.
 We have  been the source of many international collaborations. 
 
There is also an educational side: several people have told us that they were led
into mathematics through working as an editor.  Here is a typical story.

\footnotesize
Subject: Reminiscence from a young mathematician

I wanted to relay a bit of nostalgia and my heartfelt thanks.  Back in the late 1990s, 
I was a high school student in Oregon.   While I was interested in mathematics,
I had no significant mathematically creative outlet until I discovered the OEIS in the course of 
trying to invent some puzzles for myself.   
 I remember becoming a quite active contributor
through the early 2000s, and eventually at one point, an editor.  
   My experience with the OEIS, and the eventual intervention of one of 
my high school teachers, catalyzed my interest in studying mathematics, 
which I eventually did at [...] College.
I went on to a Ph.D. in algebraic geometry at the University of [...],  and am currently at [...].

I wanted to thank you for seriously engaging with an 18-year old kid, even though
I likely submitted my fair share of mathematically immature sequences.
I doubt I would have become  a mathematician without the OEIS!
\normalsize

A less-obvious use of the database is to quickly tell you how hard
a problem is. I use it myself in this way all the time.
Is the sequence ``Catalan'' or ``Collatz''?
If a sequence comes up in your own work, or when reviewing someone else's work,
it is useful to know right away if this is a well-understood sequence, 
like the Catalan numbers, or if it is one of the notoriously intractable problems 
like the Collatz or $3x+1$ problem (\seqnum{A006577}).

Finally,  the OEIS is a welcome escape when you feel the world is falling apart.
Take a look at Scott Shannon's drawings of stained glass windows
in \seqnum{A331452}; or Jonathan Wild's delicate illustrations
of the ways to draw four circles in \seqnum{A250001}; 
or \'{E}ric Angelini's ``1995'' puzzle (\seqnum{A131744}) or any of his 
``lexicographically earliest sequences" (\seqnum{A121053}, \seqnum{A307720}, and many more);
or find better solutions to the Stepping Stones Problem (\S\ref{SecStones},\seqnum{A337663}).
You can find brand new problems at any hour of the day or night
by looking at the stack of recent submissions: but beware, you may
see a problem there that will keep you awake for days. 
Or search in the database for phrases like
``It appears that ...'', or ``Conjecture: ...'',
or ``It would be nice to know more!''

\subsection{Layout of a typical entry}\label{SecEntry}

This is a good place to mention some of the features of an OEIS entry. Most of the fields (see
Figs.~\ref{FigA108a} and \ref{TableA108b}) are self-explanatory.
At the top it tells you how many matches were found to your query ($26$ in the example).  These are ranked in order of importance.

The DATA section shows the start of the sequence, usually enough terms
to fill a few lines on the screen (typically  $300$ to $500$ decimal digits).
Often one wants more terms than are shown, and
the first link in the entry will point to a plain text file with perhaps $10000$ or $20000$ terms.
That file will have a name like b001006.txt, and is called the ``b-file'' for the sequence.
Some entries also have much larger tables, giving a million or more terms.

If you  click the ``graph''  button near the top of
the reply, you will be shown two plots of the sequence,
and if you click the ``listen'' button, you can listen to the
sequence played on an instrument of your choice. The default
instrument is the grand piano, and the terms of the sequence
would be mapped to the $80$ keys by reducing the numbers mod $80$ and adding $1$.

I conclude this section with a  philosophical comment.

When you are seriously trying to analyze a sequence, and are prepared to spend any amount of time 
needed (searching for a formula or recurrence, for instance), you need all the help you can get, which is why
we provide the b-files and other data files, and why we give computer programs in so many languages. This is also
the reason we give as many references and links as possible for a sequence.  Even if the 
reference is to an ancient or obscure journal, or  one that has been accused as being ``predatory'',
we still give the reference, especially for sequences that are not well-understood.  
The same thing holds for formulas, comments, and cross-references
to other sequences.  When you are desperate, you will accept help from anywhere.
And do not forget  \href{https://oeis.org/ol.html}{``Superseeker''!}

\subsection{Arrangement of the entries}\label{SecArrange}

The entries in the database are (virtually) arranged in two different ways,
the first  essentially chronological, the second lexicographic.

The first is by their {\em absolute} identification  number, or A-number.\footnote{The
sequences in the 1973 and 1995 books were numbered N0001,... and M0001,... respectively.}
Once the collection reached a few hundred entries, I sorted them into lexicographic order 
and numbered them A1, A2, A3, $\ldots$.  
A1 gives the number of symmetry groups of order $n$, A2 is the 
famous Kolakoski sequence, and so on.  This numbering is still used today, 
only A1 has become \seqnum{A000001}, 
A2 is \seqnum{A000002}, ..., and as each new submission comes in it gets a number from the stack.
Current sequences are being issued numbers around A360000. 
Rejected A-numbers are recycled, so there are no gaps in the order.
We reached $100000$ entries in 2004, and $250000$ in 2015.
The present growth rate is about $12000$ new entries each year.

The second arrangement is a kind of lexicographic ordering.
First I describe an idealized, theoretical, lexicographic order.
Sequences of nonnegative numbers can be arranged in lexicographic (or dictionary)
order.  For example, sequences beginning $1, 2, 4, \ldots$ come before
$1, 2, 5, \ldots$, $1, 2, 4, 3, \ldots$, $1, 3, \ldots$, etc., but after $1, 2, 3, \ldots$.
Also $1, 2, 4, \ldots$ comes after
the two-term sequence $1, 2$ (because blanks precede numbers).

More formally, we compare the two sequences term-by-term, and in the first position
where they differ whichever is smaller (or blank) is the lexicographically earlier sequence.

For sequences with negative terms, we ignore the signs and sort
according to the absolute values.

Here is the actual ordering used in the OEIS.
The  sequences are arranged (virtually) into a version of
lexicographic order, according to the following rules.
First, delete all minus signs.
Then find the first term that is greater than $1$,
and discard all the terms before it. What's left determines its position
in the lexicographic order.
For example, to place $-1, 0, 1, 1, \underline{2}, 1, 17, 3, 2, 1, \ldots$
in the ordering, we would ignore the terms before the underlined $2$,
and consider the sequence as beginning  $2, 1, 17, 3, 2, 1, \ldots$.

Sequences that contain only $0$s, $1$s and $-1$s are sorted into lexicographic order
by absolute value and appear at the beginning of the ordering.
The first sequence in the database is therefore the zero sequence \seqnum{A000004}.

In this way every sequence has a unique position in the ordering.
The sequences have been sorted in this way since the 1960s.
For the first ten years the punched card entries  were physically sorted into this order.

When you look at an OEIS entry, \seqnum{A005132} say (the subject of Section~\ref{SecRecaman}), 
towards the bottom you will see two lines
like\footnote{If you don't see these, click on the A-number at the top of the entry.}

Sequence in context: A277558 A350578 A335299 * A064388 A064387 A064389

Adjacent sequences:  ~A005129 A005130 A005131 * A005133 A005134 A005135

\noindent
which tell you the three entries immediately before and after that entry in the lexicographic ordering,
and the three entries before and after it in the A-numbering. 
The asterisks represent the sequence you are looking at.
The first group can be useful if you are uncertain about a term in your sequence, 
the second in case you want to look at other sequences submitted around that time.

Today the sequences are actually  stored internally in an  SQLite database. However, the 
punched card format has been so useful that when you view a sequence, as in Fig.~\ref{TableA108b}, 
 it is still presented to you in something very like the old punched card format.


\subsection{Summary: ``A Handbook of Integer Sequences" today}\label{SecFacts}
\begin{list}{--}{\setlength{\itemsep}{0.02in}}
\item Now  {\em The On-Line Encyclopedia of Integer Sequences} or {\em OEIS}:
\href{https://oeis.org}{\tt https://oeis.org} 
\item Accurate  information about 360000 sequences.
\item Definition, formulas, references, links, programs.
View as list, table, graph, music!
\item Traffic: 1 million hits/day.
\item $30$ new entries, $50$ updates every day. 
\item Often called one of best math sites on the Web. Fingerprint file for mathematics.
\item Street creds: \href{https://oeis.org/wiki/Works_Citing_OEIS}{10000 citations}.
\item A moderated Wiki, owned by OEIS Foundation, a 501(c)(3) public charity.
\item Uses: to see if your sequence is new, to find references, formulas, programs.
\item Catalan or Collatz?   (Very easy or very hard?)
\item Source of fascinating research problems;\footnote{Look for ``Conjecture'', ``It appears that'', ``It would be nice to'', ...} low-hanging fruit from recent submissions.
\item Accessible (free, friendly).
\item Fun ($1, 2, 4, 6, 3, 9, 12, 8, 10, 5, 15, ...$?). Interesting, educational.  Escape.
\item Addictive (better than video games).
\item Has led many people into mathematics.
\item One of the most successful international collaborations, a modest contribution towards world peace.
\item Need editors.
 \end{list}

%
%
\section{Some favorite sequences}\label{SecFavorite}

I'm sometimes asked what my favorite sequence is.  This is a difficult question.
I'm tempted to reply by saying: If you were the keeper of
the only zoo in the world, how would you answer that question?
(Because that is roughly  the situation I'm in.) Would you pick one of
the exotic animals, a giraffe, a kangaroo, or a blue whale?
Or one of the essential animals, like a horse, a cow, or a duck?
If the question came from a visiting alien, of course, there is
only one possible answer: a human being.

For sequences, the essential ones are the primes, the powers of $2$,
the Catalan numbers, or (especially if the question
 came from an alien with no fingers or toes),
the counting sequence 0, 1, 2, 3, 4, ... (\seqnum{A001477}).

But here I'll mention a few that are fairly exotic.
The Recam\'{a}n and Gijswijt sequences have simple
recursive  definitions, yet are astonishingly hard to understand.

\subsection{Recam\'{a}n's sequence (A005132)}\label{SecRecaman}
This remarkable sequence has resisted analysis for over $30$ years,
even though we have computed an astronomical number of terms.
It was contributed to the database by Bernardo Recam\'{a}n Santos in 1991.

The definition is deceptively simple. The first term is $0$.
We now add or subtract $1$, then we add or subtract $2$, then add or
subtract $3$, and so on. The rule is that we always first try to
subtract, but we can only subtract if that leaves a nonnegative
number that is not yet in the sequence. Otherwise we must add.

Here is how the sequence starts. We have the initial $0$. We can't subtract $1$,
because that would give a negative number, so we add $1$ to $0$.
So the second term is $1$.
We can't subtract $2$ from $1$, so we add it, getting the third term $1+2 = 3$.
Again we can't subtract $3$, for that would give $0$, which has already appeared, so we add $3$, getting
the fourth term $3+3 = 6$.

Now we must add or subtract $4$, and this time we can subtract,
because $6-4 = 2$, and $2$ is nonnegative and a number that hasn't yet appeared.
So at this point the sequence is $0, 1, 3, 6, 2$.  Then it continues $7 (= 2+5)$,
 $13 (=7+6), 20 (=13+7), 12 (=20-8)$, and so on.  The first 16 terms are

                        $$0, 1, 3, 6, 2, 7, 13, 20, 12, 21, 11, 22, 10, 23, 9, 24, \ldots$$

When adding rather than subtracting, repeated terms are permitted ($42$ is repeated at the $24$th term).

Edmund Harriss has found an elegant way to draw the sequence as a spiral on the number line.
Start at $0$, and when we subtract $n$, draw a semicircle of diameter $n$
to the left from the last point, or to the right if we are adding $n$.
Draw the semicircles alternately below and above the
horizontal axis so as to produce a smooth spiral.

        \begin{figure}[!ht]
  \centerline{\includegraphics[ width=5in]{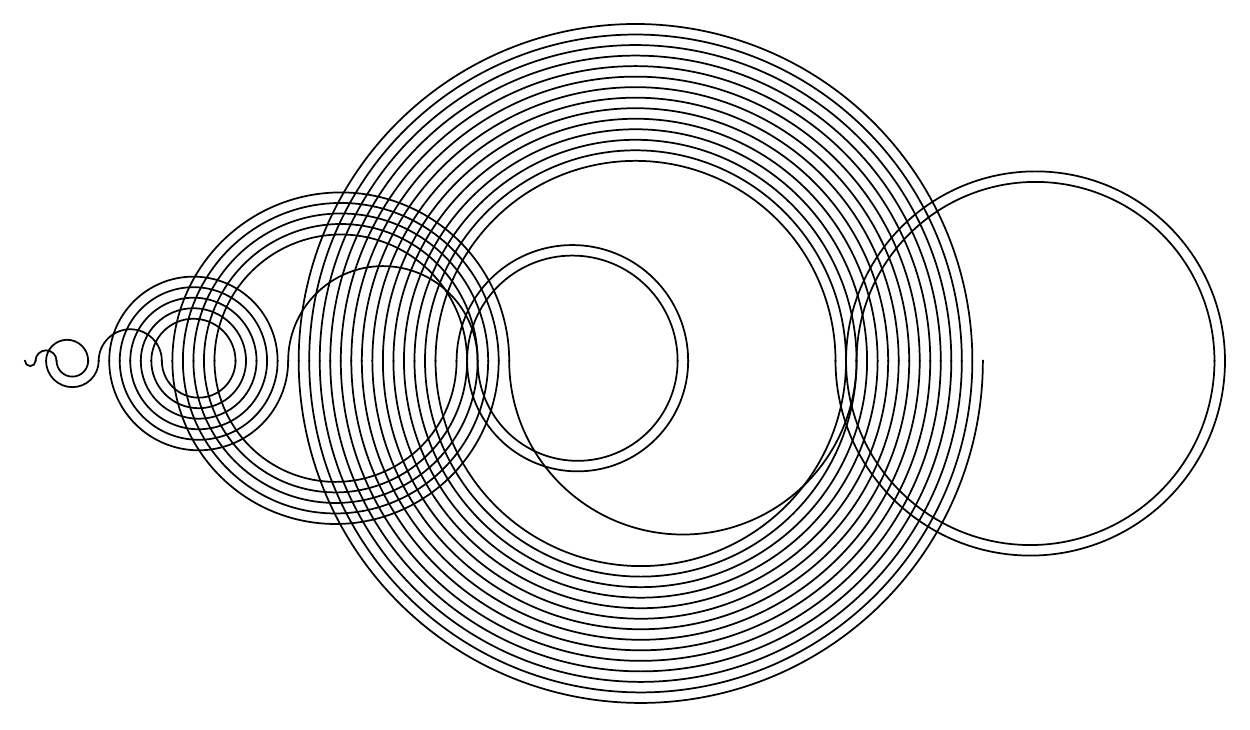}}
  \caption{Harriss's drawing of the first $64$ terms of Recam\'{a}n's sequence. (The tiny initial semicircle, 
at the extreme left, is below the axis. It has diameter $1$ and joins the points $0$ and $1$. It continues as a semicircle of diameter $2$, above the axis, joining the points $1$ and $3$.)}
        \label{Figrecaman1}
        \end{figure}

The main question about this sequence is: Does every positive number appear?  
What makes this sequence
so interesting is that certain numbers (for reasons we do not
understand) are extremely reluctant to appear.  $4$ does not appear
until $131$ steps, and $19$ takes $99734$ steps.

A group of us at AT\&T Bell Labs worked on this sequence in 2001,
and developed a way to greatly speed up the computation.
Allan Wilks used it to compute
the first $10^{15}$ terms,
and found that
$2406$ (which had been missing for a long time) finally appeared
at step $394178473633984$.

At this point the smallest missing number was $852655 = 5 \cdot 31 \cdot 5501$. Benjamin Chaffin has continued this work, and in 2018
reached $10^{230}$ terms.
However, $852655$ was still missing, and there has been no progress since then.

Thirty years ago I thought that every number would eventually appear.
Now I am not so sure. My current belief is that there are two possibilities: either there are infinitely many 
numbers that never appear, and $852655$ just happens to be the smallest of them,
but has no other special property.  A similar  phenomenon seems to occur when iterating various number-theoretic functions---see the next section.
Or, every number will eventually appear (just as presumably every one of Shakespeare's plays
will eventually appear in the expansion of $\pi$ in base $60$), although we may never
be able to extend the sequence far enough to hit $852655$.
For the latest information about this sequence (or any other sequence mentioned in this article),
consult the OEIS.

Open question: Does $852655$ appear in \seqnum{A005132}?

\subsection{Iteration of number-theoretic functions}\label{SecFunWithDigits}
Many mysterious sequences arise from the
iteration of number-theoretic functions.
A classic problem concerns the iteration of the function
$f(n) = \sigma(n) - n$, the sum of the ``aliquot parts'' of $n$ (see Guy \cite[\S B6]{UPNT}, \seqnum{A001065}).
For an initial value of $n$, what happens to the trajectory $n, f(n), f(f(n)), \ldots$?
All $n < 276$ terminate by entering a cycle (such $n$ are called  ``perfect'', ``amicable'', or ``sociable'' numbers),
or  reaching a prime, then $1$, then $0$. But it appears likely  that $n = 276$ 
and perhaps all sufficiently large even numbers,
will never terminate \cite{CGJ20}. The trajectory of $276$ is sequence \seqnum{A008892}. 
At the time of writing, the trajectory has been computed for $2145$ terms, and is still growing, term $2145$  
being a $214$-digit number \cite{FDB}. 
\seqnum{A098007} gives the number of distinct terms in the trajectory of $n$, or $-1$ if the trajectory is unbounded. The value of \seqnum{A098007}$(276)$ is unknown.

If indeed $276$ does go to infinity, it is natural to ask, how did $276$ know it
was destined to be the first immortal number under the map $f$?
The answer may be that there are infinitely many immortal numbers, and $276$ just happens to be the first.  It got lucky, that's all!  Just as $852655$ got lucky in Recam\'{a}n's problem.

A similar question, also discussed by Guy \cite[\S B41]{UPNT}, which has received much less attention, concerns
the map $g(n) = (\sigma(n) + \phi(n))/2$, 
where $\phi(n)$ is the Euler totient function \seqnum{A000010}.
The trajectory may end at $1$, a prime, or a fraction, or it may increase monotonically to infinity.
Sequence \seqnum{A292108} gives the number of steps in the trajectory, or $-1$ if the trajectory is infinite.
All numbers $n < 270$ have finite trajectories, but it appears that $270$ goes increases forever.
 The trajectory of $270$ is \seqnum{A291789}.  Andrew Booker has given a heuristic argument showing
 that almost all numbers go to infinity.  What makes $270$ the first immortal number under $g$?
Again I suspect it just got lucky!

Open questions: Does the trajectory of $276$ under $f$ increase forever? 
What about the trajectory of $270$ under $g$?

\subsection{Gijswijt's sequence (A090822)}\label{SecGijs}
For this sequence it will be helpful to remember that  chemists do not write 
$H-H-O$, they write $H_2 \,O$, they do not write 
$Al Al Al S O O O O S O O O O$, they write $Al_3 (SO_4)_2$.
We will apply a similar compression to
sequences of numbers, except that we indicate repetition 
by superscripts rather than by subscripts. 

For this problem, when we look at a sequence of numbers, we want to write
it in the form $X Y Y  \ldots Y$, or $X Y^k$, where $X$ and $Y$ are themselves
sequences of numbers, $X$ can be missing, and the exponent $k$ is {\em as large as possible}.

For example, we can  write  $1, 2, 2, 2,2$  as  $X Y^k$, where $X = 1$, $Y = 2$, and $k = 4$.
The highest $k$ we can achieve for a sequence is called its  ``curling number''.
So $1,2,2,2,2$ has curling number $4$.
Think of an animal with its head looking to the left, with a very curly tail. 
$X$ represents the head and body of the animal, and $Y^k$ represents the curls
in its tail.

Consider the sequence $3,2,4,4,2,4,4,2,4,4$. We could take 
$X = 3,2,4,4,2,4,4,2$ and $Y = 4$, getting $X Y^2$, with $k=2$,
or we could take $X = 3$, $Y = 2,4,4$, getting $X  Y^3$, with $k=3$,
which is larger.  So this sequence has curling number $3$.

Remember that  $X$ may be missing. So the sequence with a single term $99$,
say, can be written as $Y^1$ where  $Y$ is the number $99$, and it has curling number $1$.
The notion of curling number is independent  of the base in which the numbers are written.

We are now ready to define Dion Gijswijt's absolutely brilliant sequence,
which he sent to the OEIS in 2004.

The rule for finding the next term is simple: it is the curling number of
the sequence so far.  And you start with $1$. 
That's the sequence!

So let's construct it. We start with $1$, and the curling number of $1$ is $1$.
So now we have $1, 1$.  This has curling number $2$, so now we have $1, 1, 2$.
At each step we recompute the curling number, and make that the next term.

Here are the first few generations.

\begin{verbatim}
1
1 1
1 1 2
1 1 2 1
1 1 2 1 1
1 1 2 1 1 2
1 1 2 1 1 2 2 (we took Y = 1 1 2)
1 1 2 1 1 2 2 2
1 1 2 1 1 2 2 2 3
\end{verbatim}
and we have found the first $3$, at the $9$th term.
After a while, a $4$ appears at term $220$.

But Gijswijt was unable to find a $5$, and left that question open
when he submitted the sequence.
Some Bell Labs colleagues computed many millions of terms,
but no $5$ appeared.

Finally, over the course of a long weekend, Fokko van der Bult
(a fellow student of Gijswijt's in Amsterdam) and I independently
showed that there is a $5$.  In fact there are infinitely many $5$'s,
but the first one does not appear until about term $10^{10^{23}}$.
The universe would be cold long before any computer search would find it. 

In the paper we wrote about the sequence \cite{BGL07}, we also conjectured that
the first time a number $N > 4$ appears is at about term
$$2 \uparrow (2 \uparrow (3 \uparrow (4 \uparrow (5 \uparrow \ldots \uparrow (N-1))))),$$
where the up-arrows  ($\uparrow$) indicate exponentiation.
This is a tower of exponents of height $N-1$.

A very recent manuscript by a student of Gijswijt's, 
Levi van de Pol \cite{Pol22}, still under review, has extended our work, 
and may have proved the above conjecture.

I cannot resist adding a further comment about curling numbers, 
which if true shows that the Gijswijt sequence is in a sense universal.

The Curling Number Conjecture asserts that if any finite starting sequence is extended by 
the rule that the next term is the curling number of the sequence so far, 
then eventually the curling number will be $1$.

If true, this implies that if the starting sequence contains no $1$s, then the sequence eventually becomes Gijswijt's sequence \cite[Th. 23]{CLSW13}. In fact I conjecture that this is true for any starting sequence.

Open question:  Is the Curling Number Conjecture true?

\subsection{Lexicographically Earliest Sequences}\label{SecLES}

Although there is no space to discuss them in detail, let me just mention that there are
many fascinating and difficult  sequences in the OEIS whose definition has the 
form ``Lexicographically Earliest Sequence 
of distinct positive numbers with the property that ...", where now we are 
using  lexicographic in its pure sense, as defined 
in Section~\ref{SecArrange}. A favorite example is the EKG (or ECG) sequence \seqnum{A064413},
whose definition is the lexicographically earlier infinite sequence of distinct positive numbers with the property that 
each term after the first has a  nontrivial common factor with the previous term \cite{LRS02}.
Other L.E.S. examples are the Yellowstone permutation \seqnum{A098550} \cite{AHS15},
the Enots Wolley sequence \seqnum{A336957} (the name suggests the definition),
and the Binary Two-Up sequence \seqnum{A354169} \cite{DSS22}.

Open question:  Show that the terms of the Enots Wolley sequence are precisely 1, 2, and all numbers with at least two distinct prime factors.

\subsection{The Stepping Stones Problem (A337663)}\label{SecStones}

This lovely problem was invented in 2020 by two undergraduates,
Thomas Ladouceur and Jeremy Rebenstock. You have an infinite chessboard, and a handful of brown
stones, which are worth one point each. You also have an infinite number of
white stones, of values $2$, $3$, $4, \ldots$, one of each value.
Suppose you have $n$ brown stones. You start by placing them anywhere on the board.
Now you place the white stones, trying to place as many as you can.
The rules are that you can only place a white stone labeled $k$ on a square
if the values of the stones on the eight squares around it add up to $k$.
And you must place the white stones in order, first $2$, then $3$, and so on.
You stop when you cannot place the next higher-numbered white stone.
The goal is to maximize the highest value that you place. Call this $a(n)$.


\begin{figure}[!ht]
$$
\begin{array}{|c|c|c|c|c|c|c|c|}
\hline
  ~~~ & ~~~   & ~~~  &  ~~~ &  ~~~ & ~~~  &  ~~~ & ~~~  \\
\hline
 ~~~\, & ~9~ & 5 &10 &\,11\, &  ~~ &  ~~ & ~~~ \\
\hline
  &   & 4 & \setlength{\fboxrule}{2pt}\fbox{1} &   &   &   &   \\
\hline
  & 12& ~8~ & 3 & 2 &   &\,16^{*} &   \\
\hline
  &   &   &   & 6 & \setlength{\fboxrule}{2pt}\fbox{1}  &15 &   \\
\hline
  &   &   &13 & 7 &14 &   &   \\
\hline
  &   &   &   &   &   &   &   \\
\hline
\end{array}
$$
 \caption{A solution to the Stepping Stones problem for two starting stones.
The high point $a(2)=16$ here  is indicated by an asterisk, as it is in the next three tables.}
  \label{Table2.16}
\end{figure}

Say we start with $n=2$ brown stones. There are infinitely many squares where they can be placed,
 but it turns out that the best thing is to place them so they are separated  
 diagonally by a single blank square, as in Fig.~\ref{Table2.16}.
 Now we start trying to place the white stones.
 The $2$ stone has to go between the two brown (or $1$) stones, and then the $3$ 
 goes on a square adjacent to the $1$ and the $2$.
 There is now a choice for where the $4$ goes, but the choice shown
 in Fig.~\ref{Table2.16} is the best.
(After we have placed the $4$, the neighbors of the $3$ no longer add up to $3$,
but that is OK. It is only when we {\em place} the $3$ that its neighbors must add to $3$.)
Continuing in this way, we eventually reach $16$.  There is nowhere to place the $17$, so we stop.
Ladouceur and Rebenstock showed, using a computer and considering all possible arrangements,
that $16$ is the highest value that can be attained with two starting stones.
So $a(2)=16$.

This is clearly a hard problem, since the number of possibilities grows rapidly
with the number of brown stones. Only six terms of this sequence are known:
$a(1)$ through $a(6)$ are $1, 16, 28, 38, 49, 60$.
A solution for $n=4$ found by Arnauld Chevallier
 is shown in Fig.~\ref{Table4.38}.
There are lower bounds for larger values of $n$ which may turn
out to be optimal.  For  $n = 7,\ldots,10$ the current best constructions give
$71, 80, 90, 99$.
See \seqnum{A337663} for the latest information.


\begin{figure}[!ht]
$$
\begin{array}{|c|c|c|c|c|c|c|c|c|c|c|c|c|c|c|}
\hline
~~~ & ~~ & ~~ & ~~ & ~~ & ~~ & ~~ & ~~ & ~~ & ~~ & ~~ & ~~ & ~~ & ~~ & ~~~ \\
\hline
~~ & 35 & 18 & 36 &  ~~ & 23 &  ~~ & 21 &  ~~ & 32 &  ~~ &  ~~ &  ~~ &  ~~ &  ~~ \\
\hline
~~ &  ~~ & 17 &  \setlength{\fboxrule}{2pt}\fbox{1}  &  ~~ & 14 &  9 &  ~~ & 12 & 20 &  ~~ &  ~~ &  ~~ &  ~~ &  ~~ \\
\hline
~~ &  ~~ & 34 & 16 & 15 &  ~~ &  5 &  4 &  8 &  ~~ &  ~~ & 26 & 27 &  ~~ &  ~~ \\
\hline
~~ &  ~~ &  ~~ &  ~~ & 31 &  ~~ & 10 &  \setlength{\fboxrule}{2pt}\fbox{1}  &  3 & 19 & 25 &  ~~ &  \setlength{\fboxrule}{2pt}\fbox{1}  & 28 &  ~~ \\
\hline
~~ &  ~~ &  ~~ &  ~~ &  ~~ &  ~~ & 11 &  ~~ &  2 &  6 &  ~~ & 33 &  ~~ & 29 &  ~~ \\
\hline
~~ &  ~~ &  ~~ &  ~~ &  ~~ &  ~~ & 24 & 13 & 22 &  \setlength{\fboxrule}{2pt}\fbox{1}  &  7 &  ~~  &  ~~  &  ~~  &  ~~  \\
\hline
~~ &  ~~ &  ~~ &  ~~ &  ~~ &  ~~ & 37 &  ~~ &  ~~ & 30 & 38^{*} &  ~~ &  ~~ &  ~~ &  ~~ \\
\hline
~~ & ~~ & ~~ & ~~ & ~~ & ~~ & ~~ & ~~ & ~~ & ~~ & ~~ & ~~ & ~~ & ~~ & ~~ \\
\hline
\end{array}
$$
\caption{A solution to the Stepping Stones problem for four starting stones.}
\label{Table4.38}
\end{figure}

We don't know how fast $a(n)$ grows. There have been a series of upper and lower bounds, initiated by
Robert Gerbicz and Andrew Howroyd. 
The simple linear construction
shown in Fig.~\ref{Table4.18} shows that $a(n) \ge 6(n-1)$ for $n \ge 3$.


 \begin{figure}[!ht]
   $$
   \begin{array}{|c|c|c|c|c|c|c|c|c|c|c|c|}
   \hline
    ~~~  &  ~~  &  ~~  &  ~~  &  ~~  &  ~~  &  ~~  &  ~~  &  ~~  &  ~~  &  ~~  &  ~~~  \\
   \hline
    ~~  &   \setlength{\fboxrule}{2pt}\fbox{1}  &  ~~  &  ~~  &  ~~  &  ~~  &  ~~  &  ~~  &  ~~  &  ~~  &  ~~  &  ~~  \\
   \hline
    ~~  &  ~~  &  2  &  3  &  4  &  5  &  6  &  7  &  8  &  9  &  ~~  &  ~~  \\
   \hline
    ~~  &  ~~  &  ~~  &   \setlength{\fboxrule}{2pt}\fbox{1}  &  ~~  &  ~~  &  \setlength{\fboxrule}{2pt}\fbox{1}  &  ~~  &  ~~  &   \setlength{\fboxrule}{2pt}\fbox{1}  &  10  &  ~~  \\
   \hline
    ~~  &  ~~  &  18^{*}  &  17  &  16  &  15  &  14  &  13  &  12  &  11  &  ~~  &  ~~  \\
   \hline
    ~~  &  ~~  &  ~~  &  ~~  &  ~~  &  ~~  &  ~~  &  ~~  &  ~~  &  ~~  &  ~~  &  ~~  \\
   \hline
   \end{array}
   $$
    \caption{Every additional $1$ on the middle row increases the number of white stones by $6$,
    showing that $a(n) \ge 6(n-1)$ for $n \ge 3$.}\label{Table4.18}
   \end{figure}
   
By combining the constructions of Figs.~\ref{Table2.16} and \ref{Table4.18},
Menno Verhoeven  obtained $a(n) \ge 6n+3$ for $n \ge 3$ (Fig.~\ref{Table_4}). 

 \begin{figure}[!ht]
   $$
   \begin{array}{|c|c|c|c|c|c|c|c|c|c|c|}
   \hline
    ~~~  &  ~~  &  ~~~  &  ~~  &  ~~  &  ~~  &  ~~  &  ~~  &  ~~  &  ~~  &  ~~~  \\
   \hline
    ~~  &  ~~  &  ~~  &  ~~  &  ~~  &  ~~  &  ~~  &  ~~  &  25  &  ~~  &  ~~  \\
   \hline
    ~~  &  ~~  &  ~~  &  ~~  &  ~~  &  ~~  &  ~~  &  24  &  \setlength{\fboxrule}{2pt}\fbox{1}  &  26  &  ~~  \\
   \hline
    ~~  &  ~~  &  ~~  &  ~~  &  ~~  &  ~~  &  ~~  &  23  &  ~~  &  27  &  ~~  \\
   \hline
    ~~  &  ~~  &  ~~  &  ~~  &  ~~  &  ~~  &  ~~  &  22  &  ~~  &  28  &  ~~  \\
   \hline
    ~~  &  ~~  &  ~~  &  ~~  &  ~~  &  ~~  &  ~~  &  21  &  \setlength{\fboxrule}{2pt}\fbox{1} &  29  &  ~~  \\
   \hline
    ~~  &  ~~  &  ~~  &  ~~  &  ~~  &  ~~  &  ~~  &  20  &  ~~  &  30  &  ~~  \\
   \hline
    ~~  &  ~~  &  ~~  &  ~~  &  ~~  &  ~~  &  ~~  &  19  &  ~~  &  31  &  ~~  \\
   \hline
    ~~  &   9  &   5  &  10  &  11  &  ~~  &  ~~  &  18  &  \setlength{\fboxrule}{2pt}\fbox{1} &  32  &  ~~  \\
   \hline
    ~~  &  ~~  &   4  &   \setlength{\fboxrule}{2pt}\fbox{1}  &  ~~  &  ~~  &  ~~  &  17  &  ~~  &  33^{*}  &  ~~  \\
   \hline
    ~~  &  12  &   8  &   3  &   2  &  ~~  &  16  &  ~~  &  ~~  &  ~~  &  ~~  \\
   \hline
    ~~  &  ~~  &  ~~  &  ~~  &   6  &  \setlength{\fboxrule}{2pt}\fbox{1}  &  15  &  ~~  &  ~~  &  ~~  &  ~~  \\
   \hline
    ~~  &  ~~  &  ~~  &  13  &   7  & 14  &  ~~  &  ~~  &  ~~  &  ~~  &  ~~  \\
   \hline
    ~~  &  ~~  &  ~~  &  ~~  &  ~~  &  ~~  &  ~~  &  ~~  &  ~~  &  ~~  &  ~~  \\
   \hline
   \end{array}
   $$
    \caption{Combining the constructions of
    of  Figs.~\ref{Table2.16} and \ref{Table4.18}
gives $a(n) \ge 6n+3$ for $n \ge 3$.
The case $n=5$ is shown. For other values of $n$, adjust the height of the ``chimney'' on the right. }
\label{Table_4}
   \end{figure}

The best lower bound for large $n$ is due to Robert Gerbicz, who has shown by a remarkable extension of 
the construction in Figs.~\ref{Table4.18} and \ref{Table_4} that
$\varliminf_{n \to \infty} a(n)/n > 6$.
(A preliminary version of his bound gives
$a(n) > 6.0128\,n-5621$ for all  $n$,
although the exact values of the constants have not been confirmed.)
In his construction the ``chimney''  on the right of Fig.~\ref{Table_4} gets expanded into a whole trellis.

One might think that with a sufficiently clever arrangement, perhaps extending the construction in
Fig.~\ref{Table4.18} so that the path wraps around itself in a spiral, one could achieve large numbers
with only a few starting stones. But a simple counting argument due to Robert Gerbicz shows this is impossible.
The current best upper bound is due to Jonathan F. Waldmann,
who has shown that $a(n) < 79 n + C$ for some constant $C$. 
See \seqnum{A337663} for the latest information, including proofs of of the  results mentioned here.

Open question: Improve the lower and upper bounds on $a(n)$. The lower bound looks especially weak.

\subsection{Stained glass windows}\label{SecScott}

\input{InputBC42}

In 1998 Poonen and Rubinstein \cite{PoRu98}  famously determined the numbers of vertices and cells in the planar graph formed from a regular $n$-gon by joining every pair of vertices by a chord. The answers are in \seqnum{A006561} and   \seqnum{A007678}.
Lars Blomberg, Scott Shannon, and I have studied versions of this question when the regular $n$-gon
is replaced by other polygons, for instance by a square in which $n$ equally-spaced points
are placed along each side and  each pair of boundary points is joined by a chord.
We also studied rectangles, triangles, etc.
In most cases we were unable to find formulas for the numbers of vertices or cells, but we collected a lot
of data, and the graphs, when colored, often resemble stained glass windows (see \cite{BSS21}
and the illustrations in \seqnum{A331452} and 
other sequences cross-referenced there).\footnote{There is no fee for downloading  images from the OEIS,
but if you use any of them, please credit the source!}
So we consoled ourselves with the motto: if we can't solve it, make art!

The most promising case to analyze seemed to be the $n \times 2$ grid (although we did not succeed 
even there).

Open question: How many vertices and cells are there in the graph for the  $n \times 2$ grid, as illustrated for $n=4$ in Fig.~\ref{FigBC42}?  Sequences \seqnum{A331763} and \seqnum{A331766} give the first $100$ terms, 
yet even with all that data we have not found a formula.

\begin{figure}[!ht]
  \centerline{\includegraphics[ width=5in]{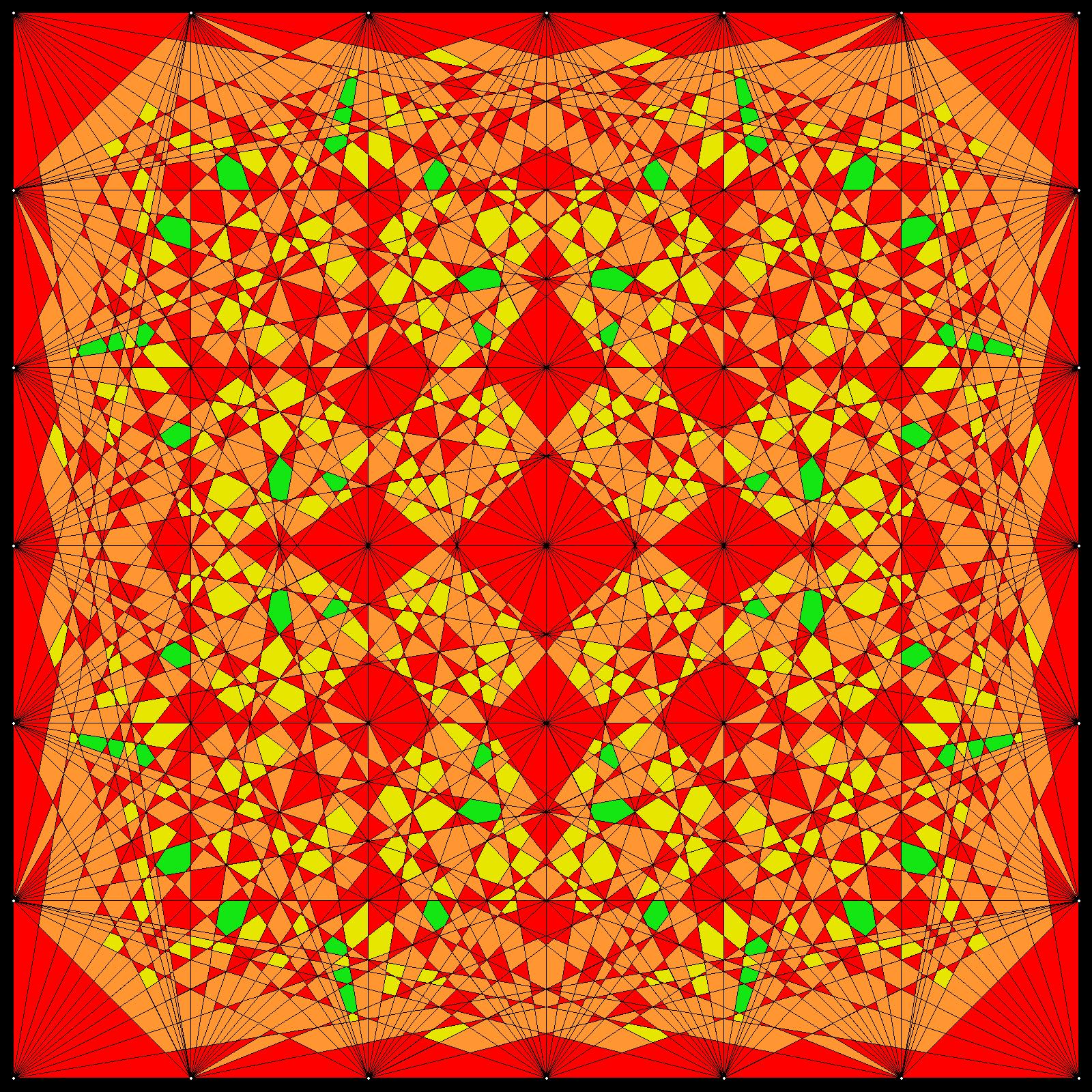}}
  \caption{A $6 \times 6$ grid with every pair
  of boundary points joined by a chord. 
  There are $4825$ vertices and $6264$ cells.} \label{FigBC66}
  \end{figure}

The case of an $n \times n$ grid seems even harder. Figure~\ref{FigBC66} shows the $6 \times 6$ graph.
Sequences \seqnum{A331449} and  \seqnum{A255011} give the numbers of vertices and cells for $n \le 42$.
\seqnum{A334699} enumerates the cells by number of sides.

\begin{figure}[!ht]
\centerline{\includegraphics[ width=5in]{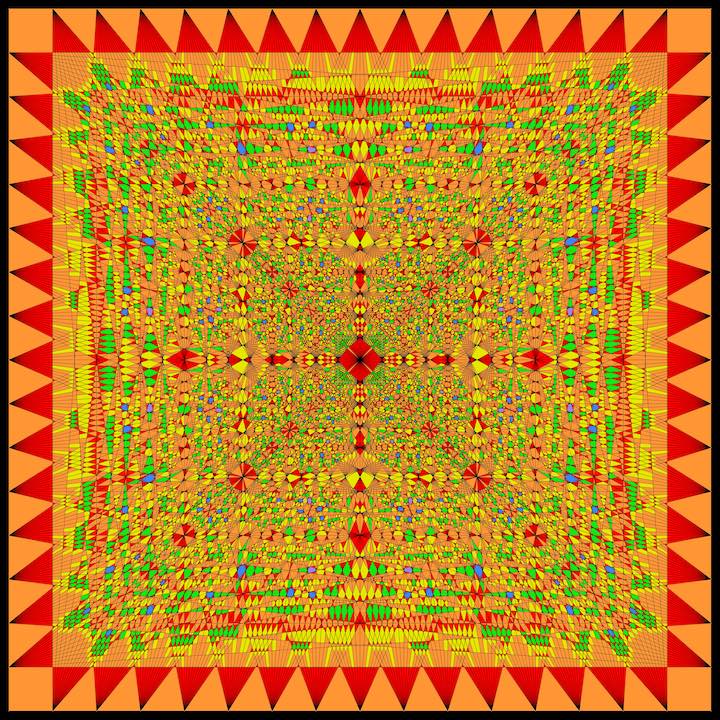}} 
 \caption{Scott Shannon's ``Magic Carpet'' graph, illustrating \seqnum{A355798}(16).}
  \label{FigCarpet}
  \end{figure}

In the summer of 2022 Scott Shannon and I considered several other families of
planar graphs. I cannot resist showing one of Shannon's graphs,
a $16 \times 16$ grid,  illustrating the $16$th term of \seqnum{A355798} (Fig.~\ref{FigCarpet}).
There are $61408$ cells. Although Shannon has  calculated $40$ terms of this sequence,
again no formula is known.

\subsection{Other sequences I would have liked to include}\label{SecSpace}

If I had had  more space I would also have discussed some very interesting sequences arising from:

-- Dissecting a square to get a regular $n$-gon (\seqnum{A110312}).

-- Gerrymandering (\seqnum{A341578}, \seqnum{A348453}, and many others).

-- In how many ways can circles overlap? (\seqnum{A250001}).

-- The Inventory sequence \seqnum{A342585}.

-- Kaprekar's junction numbers (\seqnum{A006064}, \cite{AlSl23}).

-- The kissing number problem (\seqnum{A001116}, \seqnum{A257479}).

-- The neural network problem that started it all (\seqnum{A000435}).

-- Squares in the plane (\seqnum{A051602}).

And (maybe!)  meta-sequences such as \seqnum{A051070} ($a(n)$ is the $n$th term of $A_n$) and 
\seqnum{A107357} (the $n$th term is $1 ~+~$ the $n$th term of $A_n$).

 A final comment:
there are many videos on the Internet  of talks I have given about sequences. 
There are over twenty videos that Brady Haran and I have made that have 
appeared on the Youtube Numberphile channel (and have 
been viewed over eight  million times). See for example  
\href{https://www.youtube.com/watch?v=_UtCli1SgjI}{``Terrific Toothpick Patterns"}.

\section{Acknowledgments}

I would like to thank  some good friends who have helped me and the OEIS over the years:
David L. Applegate,
William Cheswick,
Russ  Cox, 
Susanna S. Cuyler, 
Harvey P. Dale,
Ronald L. Graham, 
Richard K. Guy, 
Marc LeBrun,
John Riordan,
and
Doron Zeilberger.

There are many active volunteer editors, and it is impossible to
thank them all.  But I would like to give particular thanks to
J\"{o}rg Arndt, 
Michael S. Branicky,
Michael De Vlieger, 
Amiram Eldar,  
Charles  R. Greathouse IV, 
Maximilian  F. Hasler,  
Alois P. Heinz, 
Andrew Howroyd,
Sean A. Irvine, 
Antti Karttunen,
Michel Marcus, 
Richard J.  Mathar,
Peter Munn, 
Hugo Pfoertner,  
Kevin Ryde, 
Jon E. Schoenfield, 
R\'{e}my Sigrist,
and
Chai Wah Wu.

I also thank the members of
the Board of Trustees of the OEIS Foundation, 
past and present, for all their help,
both to me personally and to the OEIS.

Figure credits:
Figure~\ref{Fig1}(c): Clifford A. Pickover.
Figure~\ref{Figrecaman1}: Edmund Harriss.
Figures \ref{Table2.16}, \ref{Table4.38}, \ref{Table4.18}, and \ref{Table_4}
are based on communications from Thomas Ladouceur and Jeremy Rebenstock, 
Arnauld Chevallier, Skylark Xentha Murphy-Davies, and Menno Verhoeven,
respectively.
Figures~\ref{FigBC42} and \ref{FigBC66}: Lars Blomberg and Scott R. Shannon.
Figure~\ref{FigCarpet}: Scott R. Shannon.
Other figures: the author.


\bigskip
\hrule
\bigskip

\noindent 2020 Mathematics Subject Classification:  05-00, 11-00, 11Bxx, 48-00, 68-00

\bigskip
\hrule
\bigskip

\end{document}